\documentclass[12pt,leqno,oneside]{amsart}
\usepackage{mathrsfs,dsfont}
\usepackage{amsmath,amstext,amsthm,amssymb,bbm,faktor}
\usepackage{charter}
\usepackage{typearea}
\usepackage{hyperref}
\usepackage{color}
\usepackage{nccmath}
\usepackage{float}
\usepackage{tikz}
\usetikzlibrary{patterns}
\usepackage{graphicx}
\usepackage{hyperref}
\usepackage{cleveref}

\usepackage[width=6.4in,height=8.5in]
{geometry}

\def\bt{\begin{thm}}
	\def\et{\end{thm}}
\def\bl{\begin{lem}}
	\def\el{\end{lem}}
\def\bd{\begin{defn}}
	\def\ed{\end{defn}}
\def\bc{\begin{cor}}
	\def\ec{\end{cor}}
\def\bp{\begin{proof}}
	\def\ep{\end{proof}}
\def\br{\begin{rem}}
	\def\er{\end{rem}}
\def\brs{\begin{rems}}
	\def\ers{\end{rems}}

\newtheorem{thm}{Theorem}[section]
\newtheorem{prop}[thm]{Proposition}
\newtheorem{lem}[thm]{Lemma}
\newtheorem{defn}[thm]{Definition}

\newtheorem{rem}[thm]{Remark}
\newtheorem{rems}[thm]{Remarks}
\newtheorem{cor}[thm]{Corollary}

\numberwithin{equation}{section}

\title{Density Property and Composition Operators on $H(b)$ Spaces of Finitely Connected Planar Domains}

\author{Sibel \c{S}ahin} 

\address{Mathematics Department, Mimar Sinan Fine Arts University, \.{I}stanbul, Turkey}
\email{sibel.sahin@msgsu.edu.tr}

\begin{document}
	
\begin{abstract}
In this work, the density in $H(b)$ spaces of finitely connected planar domains and the boundedness of composition operators on these function spaces are studied. Density of the algebra $\mathcal{A}(D)$ is considered for both in the cases where the defining function $b$ is an extreme and non-extreme point of the unit ball of $H^\infty(D)$. In the last part boundedness of composition operators on $H(b)$ spaces is considered and as well as a generalization of the unit disk case is given, the boundedness of composition operators with generalized Blaschke symbols over finitely connected domains is characterized.

\end{abstract}
	
	\maketitle

\section{Introduction}

One of the main areas of interest in the literature of holomorphic function theory is to understand the structure of Hardy spaces, especially the Hardy-Hilbert space $H^2$. As $H^2$ is the setting of many questions and has a vast literature, the spaces inside $H^2$ are also at the center of both operator theory and holomorphic funtion theory. One of these spaces inside $H^2$ is deBranges-Rovnyak spaces and one may found a detailed study of these spaces on the unit disc case in \cite{Fri1,Fri2} and on the unit ball of $\mathbb{C}^n$ case in \cite{Sah21}. One of the questions that is often asked for a holomorphic function space is whether the algebra of holomorphic functions continuous up to boundary, namely $\mathcal{A}(D)$, is dense in the given function space or not. For deBranges-Rovnyak spaces $H(b)$, the answer of this question is highly dependent on the defining function $b\in H^\infty$ being an extreme or non-extreme point of the unit ball of $H^\infty$. For the unit disc case the answer was given affirmatively by \cite{Sar} for non-extreme $b$ in 1994. Since classical Hardy space techniques does not often work for $H(b)$ spaces, the question was open for extreme points until Aleman \& Malman gave the full picture in 2017 in \cite{AM17}. They showed that $\mathcal{A}\cap H(b)$ is dense in $H(b)$ for all $b$ in the unit ball of $H^\infty$.

In this study we are interested in deBranges-Rovnyak spaces on a finitely connected domain $D$ whose boundary is comprised of smooth Jordan curves. After giving preliminary background about Hardy spaces and $H(b)$ spaces of finitely connected domains in Section II, we will show in Section III that the density results of both Sarason and Aleman \& Malman can be extended to $H(b)$ spaces of finitely connected planar domains and although the ideas are similar the tools we use are specific to finitely connected domains. In Section IV, we will consider the composition operators on $H(b)$ spaces of finitely connected domains. In \cite{KThesis} Karaki gave a detailed study of composition operators on model spaces and as a part composition operators on $H(b)$ spaces of unit disc, also in \cite{NHK20} they investigated the boundedness of composition operators for special cases of the defining function $b$. In this last section we will first generalize their results to parabolic and hyperbolic automorphisms of unit disc and then we will turn our attention to composition operators on $H(b)$ spaces of finitely connected domains for the case where the symbol of the composition operator is a generalized Blaschke product.   
	 
\section{Preliminaries}

Let $D$ be a finitely connected domain in the complex plane with smooth Jordan boundary $\Gamma$. Then the Hardy space $H^2(D)$ is defined to be the space of holomorphic functions on $D$ such that $|f|^2$ has a harmonic majorant in $D$ and for a fixed point $z\in D$ 
$$
\|f\|_2=[u(z)]^{1/2}, ~~f\in H^2(D)
$$
where $u$ is the least harmonic majorant of $\lvert f\lvert^2$.

As in the unit disc case, Hardy space $H^2$ can also be described through boundary values [\cite{Dur}, pg:179]:

\bt
Fix a point $z_{0} \in D$, let $f(z)$ be holomorphic in $D$ and let $\bigtriangleup$ be any subdomain of $D$ containing $z_0$ whose boundary $\Gamma$ consists of a finite number of continuously differentiable Jordan curves in $D$. Then $f\in H^2(D)$ if and only if there exists a constant $M$, independent of $\bigtriangleup$, such that 
$$
\dfrac{1}{2\pi}\int_{\Gamma}\lvert f(\xi)\lvert^2\dfrac{\partial G(\xi,z_0)}{d\eta}\lvert d\xi \lvert\leq M
$$ 
where $G(\xi,z_0)$ is the Green function of $\bigtriangleup$ with pole at $z_0$.
\et

In some results we will use another class of holomorphic functions closely connected to Hardy spaces so let us define the classes $E(D)$: A holomorphic function belongs to class $E(D)$ if there is a sequence $\{\bigtriangleup_n\}$ of domains whose boundaries $\{\Gamma_n\}$ consist of a finite number of rectifiable Jordan curves such that $\bigtriangleup_n$ eventually contains each compact subset of $D$ with $\ell(\Gamma_n)<\infty$ and 
$$
\limsup_{n\rightarrow\infty} \int_{\Gamma_n}\lvert f(z)\lvert^2 \lvert dz\lvert<\infty.
$$ 
\br
Let $\omega=\psi(z)$ be any conformal mapping of $D$ onto a domain bounded by analytic Jordan curves and let $z=\varphi(\omega)$ be the inverse mapping then by [\cite{Dur}, pg:183] if $\arg\{\varphi '(\omega)\}$ is single valued in $G$, we have $f\in E^2(D)$ if and only if $f(\varphi(\omega))[\varphi'(\omega)]^{1/2}\in H^2(G)$ with $G=\psi(D)$. 
\er

Throughout this paper our main focus will be on deBranges-Rovnyak spaces $H(b)$ of a finitely connected planar domain $D$, however for the sake of completeness we will first mention the model spaces: If $A$ is a bounded operator from a Hilbert space $H_1$ into another Hilbert space $H_2$ then the space $\mathcal{M}(A)$ is defined to be the range of $A$ with the Hilbert space structure that makes $A$ a coisometry from $H_1$ onto $\mathcal{M}(A)$ [For details of these "Hilbert spaces inside Hilbert spaces", see \cite{Sar}]. Now we can define $H(b)$ spaces:

\bd
Let $b\in H^\infty(D)$ be a non-constant holomorphic function on $D$ with $\|b\|_{\infty}\leq 1$. The holomorphic function space $H(b)$ of $D$ is then defined by the inner product
$$
\|(I-T_{b}T_{\overline{b}})^{^1/2}f\|_b=\|f\|_2
$$
where $f\in (H^2(D)\ominus Ker(I-T_{b}T_{\overline{b}})^{^1/2})$ and for any $\psi\in L^2(\Gamma)$, $T_\psi:H^2(D)\rightarrow H^2(D)$ is the Toeplitz operator with symbol $\psi$.
\ed

\brs
\begin{itemize}
	\item [1)] For details of Toeplitz operators on multiply connected domains one may check \cite{Ab74}.
	\item [2)] For a Toeplitz operator $T_\psi$, the space $\mathcal{M}(T_\psi)$ is denoted shortly by $\mathcal{M}(\psi)$.
\end{itemize}
\ers

\section{Extreme and Non-Extreme Points of $H^\infty$ and Density in $H(b)$}
	
The structure of $H(b)$ spaces changes dramatically whether the defining function is an extreme point of the unit ball of the bounded holomorphic functions or not. Hence we will first consider the classification of the extreme points in the case of finitely connected domains:

\bt(\cite{GV68}, Theorem 3): A function $f\in H^\infty(D)$ is an extreme point of the unit ball of $H^\infty$ if and only if $\|f\|_{\infty}=1$ and $\int_{\Gamma}\log (1-\lvert f\lvert)d\mu=-\infty$ where $\mu$ is the harmonic measure on $\Gamma$.
\et

\underline{$\bullet$ Part 1: Density in the case $b$ is non-extreme:}

From Theorem 3.1 we understand that $b$ is a non-extreme point of the unit ball of $H^\infty(D)$ if and only if $\log (1-\lvert b\lvert^2)\in L^1(\Gamma)$ and from this equivalence we get that if $b$ is non-extreme then there is a function $a$ such that $a$ is non-vanishing in $D$ and $\lvert a\lvert^2+\lvert b\lvert^2=1$ a.e. on $\Gamma$. The function $a$ is uniquely determined up to a modulus automorphism [\cite{Ab74}, Theorem 1.12]. Then $T_{\overline{a}}T_a=1-T_{\overline{b}}T_b$ and we have $H(\overline{b})=\mathcal{M}(\overline{a})$ and using Douglas' Criterion [\cite{Sar}, I-5(i)] we see that $\mathcal{M}(a)\subset H(b)$ contractively.

Now in this setup
\begin{equation}\label[equation]{dagger}
	h\in H(b)\Leftrightarrow T_{\overline{b}}h\in H(\overline{b})\Leftrightarrow T_{\overline{b}}h\in \mathcal{M}(\overline{a})
\end{equation}
and then using the notation of Sarason, there exists a function $h^+\in H^2(D)$ with $T_{\overline{b}}h=T_{\overline{a}}h^+$ and since for any contraction $T$ inside a Hilbert space $H$ we have $H=H(T)+\mathcal{M}(T)$ for every $f_1, f_2\in H(b)$ one has
\begin{equation}\label[equation]{StarEq}
	\langle f_1,f_2\rangle_b=\langle f_1,f_2\rangle_2+\langle f_{1}^{+}, f_{2}^{+}\rangle_2.
\end{equation}

For the sake of completeness, let us also give the following easy calculation [\cite{Sar}, IV-I]:

If $h\in H(b)$ and $\varphi$ is a bounded holomorphic function then,
$$
T_{\overline{b}}T_{\overline{\varphi}}h=T_{\overline{\varphi}}T_{\overline{b}}h=T_{\overline{\varphi}}T_{\overline{a}}h^+=T_{\overline{a}}T_{\overline{\varphi}}h^+
$$	
which gives us 	
\begin{equation}\label[equation]{DoubleStarEq}
	(T_{\overline{\varphi}}h)^+=T_{\overline{\varphi}}h^+.
\end{equation}

\bl
If $b$ is a non-extreme point of $H^\infty(D)$ then polynomials are densely contained in $\mathcal{M}(\overline{a})$.
\el
\bp
First we will show that for any polynomial $p$ of degree $n$, $T_{\overline{a}}p$ is also a polynomial of degree $n$, for any $k>n$ we have
$$
\langle T_{\overline{a}}p, z^k\rangle_2=\langle p, az^k\rangle_2=0
$$ 
and hence $T_{\overline{a}}p$ is a polynomial of same degree. Since polynomials are dense in $H^2(D)$ the result follows.
\ep

Now we can state the density result for $H(b)$ spaces in the case that $b$ is a non-extreme point of $H^\infty(D)$:
\bt
If $b$ is a non-extreme point of the unit ball of $H^\infty(D)$ then polynomials are dense in $H(b)$. 
\et 
\bp
Let $h$ be a function in $H(b)$ such that $h\perp \mathcal{M}(\overline{a})$ with respect to the inner product of $H(b)$. By [\cite{Sar}, II-7], the backward shift operator $S^*$ acts as a contraction on $H(b)$ and in fact $H(b)$ is invariant under the action of $S^*$ hence $h$ is orthogonal to $T_{\overline{a}}S^{*^{n}}h$ for all $n\geq 0$. Now by \ref{DoubleStarEq} we see that 
$$
(T_{\overline{a}}S^{*^{n}}h)^+=T_{\overline{a}}S^{*^{n}}h^+
$$
and 
$$
\langle h, T_{\overline{a}}S^{*^{n}}h\rangle_b=\langle h, T_{\overline{a}}S^{*^{n}}h\rangle+\langle h^+,T_{\overline{a}}S^{*^{n}}h^+\rangle
$$
by \ref{StarEq}. Therefore we have 
$$
\dfrac{1}{2\pi}\int_{\Gamma}a(z)(\lvert h(z)\lvert^2+\lvert h^+(z)\lvert^2)z^ndz=0.
$$
Then by [\cite{Dur}, pg:184] $a(\lvert h\lvert^2+\lvert h^+\lvert^2)\in E^1(D)$ and since $a$ is an outer function [\cite{Roy}, pg:155,(3)], $(\lvert h\lvert^2+\lvert h^+\lvert^2)\in E^1(D)$ but since this is a real valued function belonging to $E^1(D)$ we must have $h=0$. 
\ep

\underline{$\bullet$ Part 2: Density in the case $b$ is extreme:}

In this part we will extend the density result to the case when the defining function $b$ is an extreme point of the unit ball of $H^\infty (D)$. For this we will use the techniques introduced by Aleman \& Malman (\cite{AM17}) but since throughout this study we are in finitely connected domains we need to develop new tools, specifically we will use the tools constructed by Royden and Khavinson for the holomorphic function spaces defined on aforementioned type of domains:

\begin{prop}
	Let $b$ an extreme point of $H^\infty(D)$ and let $E$ be the set of points $\xi\in \Gamma$ where $\lvert b(\xi)\lvert<1$. Then for any $f\in H(b)$, there exists a unique $g\in L^2(E)$ such that 
	$$
	P_+\overline{b}f=-P_+\sqrt{1-\lvert b\lvert^2}g
	$$
	and the map $J:H(b)\rightarrow H^2(D)\oplus L^2(E)$, $Jf=(f,g)$ is an isometry with $J(H(b))^{\perp}=\{(bh,\sqrt{1-\lvert b\lvert^2}h):~h\in H^2(D)\}$.
\end{prop}  
\bp
The argument used in the proof is totally functional analytic and is verbatim to the unit disc case given in \cite{AM17} so we do not repeat in here.
\ep

Apart from the above proposition another crucial argument used is The Khintchin-Ostrowski theorem which is given as follows for the more general domains:

\bt[\cite{HJ}, pg:280]: Let $\xi_0\in D$ and $\omega_{\xi_0}$ be the harmonic measure of $D$ with respect to $\xi_0$. Let $\{ F_n\}$ be a sequence of holomorphic functions in $D$ such that 
\begin{itemize}
	\item [(i)] $\sup_n\int_{\Gamma}\log^+\lvert F_n\lvert d\omega_{\xi_0}<\infty$
	\item [(ii)] On some set $E$ of positive harmonic measure, the sequence $F_n$ converges in measure to a function $\Phi$. 
	
Then $F_n$ converges uniformly on compact subsets of $D$ to an $F$ so that $F=\Phi$ a.e. on $E$. 
\end{itemize}  

\et

Now let $N_+(D)$ be the Nevanlinna class of holomorphic functions od $D$ then we have the following inner-outer factorization result of Khavinson [\cite{Khav83}, Lemma 4.1]:

\bl
Let $f\in N_+(D)$ and $B_0(z)$ be the Blaschke product corresponding to the zeros of $f$ then 
$$
\mathcal{F}=\faktor{f}{B_0}\in N_+(D).
$$
\el 

This lemma is actually the key point in the proof of the following result which is just a generalization of a similar lemma of \cite{AM17} for the finitely connected domains (See \cite{AM17}, Lemma 4 for the details of the proof):

\bl
Let $\mathcal{A}$ be the algebra of holomorphic functions on $D$ that are continuous on $\Gamma$ and $\mathcal{C}$ be the space of Cauchy transforms of finite measures on $\Gamma$ that can be identified with the dual of $\mathcal{A}$. Suppose that $E=\{ \xi\in \Gamma:~\lvert b(\xi)\lvert<1\}$, $B=A\oplus L^2(E)$, $B'=C\oplus L^2(E)$. Then
\begin{equation}\label[equation]{AM17L4}
S=\{(\mathcal{C}\mu, h):~ \faktor{\mathcal{C}\mu}{b}\in N_+(D),~\faktor{\mathcal{C}\mu}{b}=\faktor{h}{\sqrt{1-\lvert b\lvert^2}}~a.e~on~E \}
\end{equation}
is weak-* closed in $B'$.  

\el

Having all the corresponding pieces for finitely connected domains now we can state the density property in this setting:

\bt
$\mathcal{A}\cap H(b)$ is dense in $H(b)$.
\et
\bp
Let $J:H(b)\rightarrow H^2(D)\oplus L^2(E)$ be the isometry defined in Proposition 3.4 then one obtains $J(\mathcal{A}\cap H(b))=\bigcap_{h\in H^2}Ker \ell_h$ where $\ell_h$ are the functionals corresponding to $B'$ as
$$
\ell_h=(hb, h\sqrt{1-\lvert b\lvert^2}).
$$ 
As it is pointed out in [\cite{AM17}, Theorem 1], by Hahn-Banahch theorem the weak-* closure of the set of $\ell_h$ is $J(\mathcal{A}\cap H(b))^{\perp}$ and since for $h\in H^2$, $\ell_h\in S$ given in \ref{AM17L4} we get $J(\mathcal{A}\cap H(b))^{\perp}\subset S$. Therefore, for $f\in H(b)$ and $f\in (\mathcal{A}\cap H(b))^{\perp}$ then $Jf\in S$ i.e. $Jf=(hb, h\sqrt{1-|b|^2})$ for some $h\in H^2$. So by Preposition 3.4, $Jf\in J(H(b))^\perp$ and $Jf=0$.
\ep

\section{Boundedness of Composition Operators on $H(b)$ for Special Defining Functions $b$}

In this part we will show the boundedness of composition operators on some $H(b)$ spaces for some special defining functions $b$. First of all, for non-extreme $b$, we saw in Section 3- Part 1 that there is an outer function $a$ such that $|a|^2+|b|^2=1$ a.e. on $\Gamma$. Then $(a,b)$ is called a pair or a Pythagorean pair. 

Now we will try to characterize the structure of $H(b)$ spaces for some rational, non-extreme $b$, however first we need the following result of Abrahamse,\cite{Ab79}, which is the analogue of F\'{e}jer- Riesz theorem for finitely connected domains:
\bt(\cite{Ab79}, Lemma 5): Let $\omega$ be an invertible function in $L\infty(\mu)$. If $f$ is in $(\omega H^2)^\perp\cap L^\infty(\mu)$ then $f=gh$ with $g\in H^2(D)$, $h\in (\omega H^2(D))^\perp$ and
$$
|f|=|g|^2=|h|^2.
$$
\et  
\br
In this inner-outer factorization, when $f$ is a polynomial the outer function $g$ is the polynomial in the classical F\'{e}jer-Riesz theorem.
\er
\begin{prop}
If $b$ is non-extreme, rational defining function then its Pythagorean mate $a$ is also rational.
\end{prop}
\bp
Let $b\in H^\infty(D)$ be a non-extreme point in the unit ball and be rational, then $b=\dfrac{q_1}{q_2}$ where $q_1$ and $q_2$ are polynomials with $q_2$ is non-vanishing in $D$. Since $\|b\|_\infty\leq 1$ we have $1-\lvert b\lvert^2\geq 0$ for all $\xi\in\Gamma$ and $\lvert q_2\lvert^2-\lvert q_1\lvert^2$ is a non-negative polynomial. Therefore by previous analogue of F\'{e}jer-Riesz theorem, we have a polynomial $q$ so that $\lvert q_2\lvert^2-\lvert q_1\lvert^2=\lvert q\lvert^2$. Now let $a=\dfrac{q}{q_2}$ then $a$ is rational and $\lvert a\lvert^2=1-\lvert b\lvert^2$ hence $(a,b)$ is a rational Pythagorean pair. 
\ep

Now using rational Pythagorean pairs and following the idea given in \cite{CR13}, we will characterize the structure of $H(b)$ spaces through this pair:

\bt
Let $(a,b)$ be a rational Pythagorean pair. Suppose that $\lambda_1,\dots,\lambda_n$ are the zeros of $a$ on $\Gamma$. Then the space $H(b)$ has the following structure 
$$
H(b)=\left(\left(\prod_{j=1}^{n}(z-\lambda_j)\right)H^2(D)\right)\oplus \mathbb{P}_{n-1}
$$
where $\mathbb{P}_{n-1}$ is the space of polynomials of degree at most $n-1$.
\et
\bp
By Theorem 4.1, $\lvert z^n\lvert^2+\left\lvert \prod_{j=1}^{n}(\lvert z\lvert^2-\overline{\lambda_j)}z \right\lvert^2=\lvert p(z)\lvert^2$ on $\Gamma$ for some non vanishing polynomial $p$ on $\overline{D}$. Then using the construction in \cite{CR13} set $\tilde{b}=\dfrac{z^n}{p(z)}$, $\tilde{a}=\displaystyle{\dfrac{\prod_{j=1}^{n}(\lvert z\lvert^2-\overline{\lambda_j)}z}{p(z)}}$. Clearly, $(\tilde{a},\tilde{b})$ are rational Pythagorean mates. Now since $a$ and $\tilde{a}$ have the same number of zeros on $\Gamma$, $|a|\leq C |\tilde{a}|$ a.e.on $\Gamma$ for some $C>0$ and since $\tilde{a}$ and $b$ have no common roots in $D$ by Corona theorem on finitely connected domains [\cite{MA}, pg:137] we know that there exist $f_1,f_2$ in $\mathcal{A}(D)$ such that $1=\tilde{a}f_1+bf_2$ hence by an anology to [\cite{BK87} Corollary 4, \cite{CR13} Theorem 2.3] we have $H(b)=H(\tilde{b})$. From the characterization \ref{dagger}, $f\in H(\tilde{b})$ if and only if $T_{\tilde{b}}f\in T_{\tilde{a}}H^2(D)$ and for $g\in H^2(D)$, 
$$
T_{\tilde{b}}f=T_{\tilde{a}}g\Leftrightarrow P_+(\overline{\tilde{b}}f-\overline{\tilde{a}}g)=0
$$ 
$$
\Leftrightarrow \dfrac{\overline{z^n}}{\overline{p}}f-\dfrac{\overline{\prod_{j=1}^{n}(\lvert z\lvert^2-\overline{\lambda_j)}z}}{\overline{p}}g\in (H^2(D))^{\perp}\Leftrightarrow \dfrac{\overline{z^n}}{\overline{p}}(f-\prod_{j=1}^{n}(z-\lambda_j)g)\in (H^2(D))^{\perp}
$$
$$
\Leftrightarrow f-\prod_{j=1}^{n}(z-\lambda_j)g \in (z^nH^2(D))^{\perp}
$$
Then by [\cite{Dur}, pg:182], $f-\prod_{j=1}^{n}(z-\lambda_j)g \in (z^nE(D))^{\perp}$ which by [\cite{Dur}, pg:184, eg:7] gives that $f-\prod_{j=1}^{n}(z-\lambda_j)g\in \mathbb{P}_{n-1}$. 
\ep

Using this characterization, we will first give the generalization of the main results given in [\cite{NHK20}, Theorems 3.1-3.2] for the composition operators of the unit disc case:
\bt
\begin{itemize}
	\item [(i)] Let $\varphi$ be a parabolic automorphism of the unit disc with the fixed point $\xi\in \mathbb{T}$. Then 
	$$
	C_\varphi:H(b)\rightarrow H(b)
	$$
	is bounded for $H(b)=(z-\xi)H^2\oplus \mathbb{C}$.
	\item [(ii)] Let $\varphi$ be a hyperbolic automorphism of the unit disc with fixed  points $\xi_1,\xi_2\in \mathbb{T}$. Then 
	$$
	C_\varphi:H(b)\rightarrow H(b)
	$$
	is bounded for $H(b)=(z-\xi_1)(z-\xi_2)H^2\oplus \mathbb{P}_1$.
\end{itemize}
\et
\bp
\begin{itemize}
	\item [(i)] The proof of this part is actually given in [\cite{NHK20}, Theorems 3.1-3.2, pg:4-5].
	\item [(ii)] Let $f\in H(b)$ for $H(b)=(z-\xi_1)(z-\xi_2)H^2\oplus \mathbb{P}_1$, then for some $g\in H^2(\mathbb{D})$ and $a,b\in \mathbb{C}$ we have
	$$
	f=H(b)=(z-\xi_1)(z-\xi_2)g+(az+\tilde{b}).
	$$
	Suppose that $\varphi$ is a hyperbolic automorphism of the unit disc with fixed points $\xi_1,\xi_2\in\mathbb{T}$ i.e $\varphi(\xi_1)=\xi_1$ and $\varphi(\xi_2)=\xi_2$. Now,
	$$
	f\circ\varphi=(\varphi(z)-\xi_1)(\varphi(z)-\xi_2)g(\varphi(z))+a(\varphi(z))+\tilde{b}
	$$
	$$
	=(z-\xi_1)(z-\xi_2)h_1h_2g(\varphi(z))+a(\varphi(z)-z)+az+\tilde{b}
	$$
	$$
	=(z-\xi_1)(z-\xi_2)h_1h_2g(\varphi(z))+ah_3(z-\xi_1)(z-\xi_2)+az+\tilde{b}
	$$
	$$
	=(z-\xi_1)(z-\xi_2)\underbrace{[h_1h_2g(\varphi(z))+ah_3]}_{\in H^2(D)}+(az+\tilde{b})\in H(b)
	$$
	where $h_1,h_2,h_3\in H^\infty(\mathbb{D})$ and since for each $\varphi:\mathbb{D}\rightarrow \mathbb{D}$, $C_\varphi$ is a bounded operator on $H^2(\mathbb{D})$, $g\circ\varphi\in H^2(\mathbb{D})$. Hence $C_\varphi$ is bounded on $H(b)$.
\end{itemize}
\ep

For the rest of this section we will focus on composition operators generated by finite Blaschke products over $H(b)$ spaces of finitely connected domain $D$. We will follow Royden's definition for generalized Blaschke products:
\bd
A bounded analytic function $\Phi$ on $D$ is said to be a generalized Blaschke product if 
$$
\log |\Phi|=\sum_{\nu}g(z,a_\nu)+h(z)
$$
where $g(z,\xi)$ is the Green function of $D$ and $h$ is a harmonic measure.
\ed 

Now we will state the last main result of this work:
\bt
Let $(a,b)$ and $(\tilde{a},\tilde{b})$ be Pythagorean rational pairs. Suppose that the finite, generalized Blaschke product $B$ on $D$ satisfies the following condition:\\
$B(\tilde{\lambda_j})=\lambda_j$, for all $j$, where $\tilde{\lambda_j}$'s are zeros of $\tilde{a}$ and $\lambda_j$'s are zeros of $a$.

Then $C_B:H(b)\rightarrow H(\tilde{b})$ is bounded.
\et
\bp
Let $f\in H(b)$ then $f(z)=\left(\prod_{j=1}^{n}(z-\lambda_j)\right)g+p_{n-1}(z)$ for some $g\in H^2(D)$ and $p_{n-1}\in \mathbb{P}_{n-1}$ where $\lambda_j$'s are roots of $a$. Suppose $B$ is a finite, generalized Blaschke product on $D$ such that $B(\tilde{\lambda_j})=\lambda_j$ where $\tilde{\lambda_j}$'s are roots of $\tilde{a}$. Then
\begin{equation}\label[equation]{compopdoublestar}
	f(B(z))=\left(\prod_{j=1}^{n}(B(z)-\lambda_j)\right)g(B(z))+p_{n-1}(B(z))
\end{equation}
and let us now consider \ref{compopdoublestar} in two parts, namely, $f(B(z))=f_1(z)+f_2(z)$. First of all $f_1(z)$ is a holomorphic function with zeros at $\lambda_j$'s  and also for any finite Blaschke product $B$ we have $g\circ B\in H^2(D)$ so actually $f_1(z)=\left(\prod_{j=1}^{n}(z-\tilde{\lambda_j})\right)h(z)$ for some $h\in H^2(D)$. For the second part,
$$
f_2(z)=a_{n-1}(B(z))^{n-1}+\dots+a_1B(z)+a_0
$$ 
and if we reformulate this using the open form of a finite Blaschke product (see \cite{Khav83}) and rewriting everything by adding and subtracting $\lambda_j$'s then we will actually obtain
$$
f_2(z)=\left(\prod_{j=1}^{n}(z-\tilde{\lambda_j})\right)k(z)+q_{n-1}(z)
$$
where $k(z)\in H^2(D)$ and $q_{n-1}\in \mathbb{P}_{n-1}$. Hence by the characterization theorem we have $f(B(z))\in H(\tilde{b})$ and $C_B$ is a bounded composition operator.
\ep

\end{document}